\documentclass[11pt,reqno]{amsart}
\usepackage{hyperref}
\usepackage{amsmath,amsfonts,amssymb}
\usepackage[latin1]{inputenc}
\usepackage{verbatim}
\usepackage{color}

\newtheorem{theorem}{Theorem}[section]
\newtheorem{lemma}[theorem]{Lemma}
\newtheorem{prop}[theorem]{Proposition}

\newtheorem{rmk}[theorem]{Remark}
\newtheorem{defn}[theorem]{Definition}

\newcounter{defn}

\newcommand{\sect}{\vspace{3mm} \setcounter{equation}{0} \setcounter{defn}{0} \section}

\newcommand{\w}[1]{\langle {#1} \rangle}

\newcommand{\pf}{\noindent {\bf Proof. \hspace{2mm}}}
\newcommand{\ef}{ \hfill $ \Box $ \vskip 3mm}
\newcommand{\be}{\begin{equation}}
\newcommand{\ee}{\end{equation}}
\newcommand{\bea}{\begin{eqnarray}}
\newcommand{\eea}{\end{eqnarray}}
\newcommand{\beas}{\begin{eqnarray*}}
\newcommand{\eeas}{\end{eqnarray*}}

\newcommand{\f}{\frac}
\newcommand{\na}{\nabla}

\newcommand{\bC}{{\mathbb C}}
\newcommand{\bR}{{\mathbb R}}
\newcommand{\bN}{{\mathbb N}}
\newcommand{\bZ}{{\mathbb Z}}

\newcommand{\vE}{{\mathcal E}}
\newcommand{\vF}{{\mathcal F}}
\newcommand{\vG}{{\mathcal G}}

\newcommand{\vH}{{\mathcal H}}

\newcommand{\vL}{{\mathcal L}}

\newcommand{\vT}{{\mathcal T}}

\newcommand{\vS}{{\mathcal S}}

\def\f{\frac}
\def\Dl{\Delta}
\def\na{\nabla}

\def\Dl{\Delta}

\def\i{\infty}

\begin{document}

\title[The Kramers-Fokker-Planck operator]{Spectral properties of  the Kramers-Fokker-Planck operator with a long-range potential}

\author{Xue Ping WANG}

\date{\today}

\address{Laboratoire de Mathématiques Jean Leray\\
UMR CNRS 6629\\
Nantes Universit\'{e} \\
44322 Nantes Cedex 3  France \\
E-mail: xue-ping.wang@univ-nantes.fr}


\begin{abstract}
We study real resonances and embedded eigenvalues of the Kramers--Fokker--Planck operator with a long-range potential. We prove that thresholds are only possible accumulation points of eigenvalues and that  the limiting absorption principle holds  true for energies outside an exceptional set. We also prove that the eigenfunctions associated with discrete eigenvalues decay exponentially  and those associated with embedded non-threshold  ones decay polynomially.
\end{abstract}

\maketitle

\sect{Introduction}

In this work, we consider the Kramers--Fokker--Planck (KFP, for short) operator, also called  the Kramers operator by physicists (\cite{risc}), given by
\be\label{operator}
 P=-\Delta_v+\f{1}{4}|v|^2-\f{n}{2} + v\cdot\na_x-\na V(x)\cdot\na_v, \quad (x,v) \in \bR^{2n}.
  \ee
  $V(x)$ is supposed to be  a  real-valued $C^1$ potential satisfying, for some constants $\rho>-1, C>0$ the estimate
\be \label{ass1}
|V(x)| + \w{x}|\nabla V(x)| \le C\w{x}^{-\rho},  \quad x\in \bR^n.
\ee
Here $\w{x} = (1 + |x|^2)^{1/2}$ and $n \ge 1$. In most part of this work, we assume $\rho>0$, i.e., the potential is long-range, unless in the analysis of decay of eigenfunctions associated with the discrete spectrum, the condition $\rho >-1$ is sufficient.

Let $P_0$ be the free KFP operator:
\[
P_0 = -\Delta_v+\f{1}{4}|v|^2-\f{n}{2} + v\cdot\na_x, \quad (x,v) \in \bR^{2n}.
\]
$P_0$ defined on  $\mathcal{C}_0^{\infty}(\mathbb{R}^{2n})\;$ is essentially maximal accretive  (see Proposition 5.5, \cite{hln}). We still denote by $P_0$ its closed extension with maximal domain $D(P_0)=\{f \in L^2(\bR^{2n}); P_0 f \in L^2(\bR^{2n})\}$.  One has: $\sigma(P_0) =[0, \infty[$. The following subellitptic estimate holds true for $P_0$ (\cite{hln}):
\be \label{subelliptic}
\|\Delta_v f\| +\|v^2 f\| + \|\w{D_x}^{\f 2 3}f\| \le C (\|P_0 f\| + \|f\|) \quad f \in D(P_0).
\ee
For $r\ge 0,s \in \bR$, we introduce accordingly the weighted Sobolev space:
\[
\vH^{r,s} = \{u \in \vS'(\bR^{2n}); (1 + \w{D_v}^2 + |v|^2 + \w{D_x}^{\f 2 3})^{\f r 2 } \w{x}^s u \in L^2\}.
\]
For $r <0$ and $s \in \bR$, $\vH^{r,s}$ is taken to be the dual space of $H^{-r,-s}$.  Let $\vL(r,s; r',s')$ be the space of bounded linear operators from $\vH^{r,s}$ to $\vH^{r',s'}$.
To simplify notation we denote  $\vH^s =\vH^{0,s}$, $\vH = \vH^0$ and $\vL(s,s') = \vL(0,s; 0,s')$.
\\

We write the full KFP operator $P$ as
\[
P= P_0 + W \quad \mbox{ with }W = -\na_x V(x) \cdot \na_v.
\]
Under the condition (\ref{ass1}) with $\rho>-1$, $\nabla_xV(x)$ tends to $0$ as $|x| \to \infty$.  (\ref{subelliptic}) implies that $W$ is relatively compact with respect to $P_0$. Therefore $P$ defined on $D(P) =D(P_0)$ is closed and maximally accretive. Let $\sigma (P)$ (resp., $\sigma_d(P)$, $\sigma_p(P)$) denote the spectrum (resp., the discrete spectrum, the point spectrum) of $P$. $\sigma_d(P)$ is the set of isolate eigenvalues with finite (algebraic) multiplicity. Let
$\sigma_{\rm ess}(P) =\sigma(P) \setminus \sigma_d(P)$ be the essential spectrum of $P$. When $\rho >-1$, zero may be an embedded eigenvalue and under some additional conditions, the return to equilibrium is proven in \cite{lz} with  subexpoential convergence rate.
\\

 The general theory on relatively compact perturbation (\cite{d}) for closed non-self-adjoint operators affirms that
\be
\sigma_{\rm ess} (P) = \sigma (P_0)=[0, \infty[.
\ee
and $\sigma(P) \setminus [0, \infty[ =\sigma_d(P)$ and that the discrete spectrum $\sigma_d(P)$ may accumulate towards any point in $[0, \infty[$. One of the goals of this work is to improve this statement. \\

The set of  thresholds, $\vT$,  of the KFP operators is composed of  eigenvalues of the harmonic oscillator $-\Delta_v+\f{1}{4}|v|^2-\f{n}{2}$ which are non-negative integers:
\[
\vT = \bN.
\]
Motivated by spectral analysis of non-self-adjoint Schr\"odinger operators, we introduce real (non-threshold) resonances of the KFP operator using the boundary values of the free resolvent $R_0(z) =(P_0-z)^{-1}$
\[
R_0(\lambda \pm i0) = \lim_{z \to \lambda, \pm \Im z >0} R_0(z), \quad \mbox{ in } \quad \vL( s, -s)
\]
for $s> \f 1 2$. (See Proposition \ref{prop2.1} below.) The notion of threshold resonance is dimension-dependent and requires faster decay of potentials. For the KFP operator with a fastly decaying potential, see \cite{nw,w3} when $n=1$ and $3$.\\

\begin{defn}\label{defn1.1} Assume condition (\ref{ass1}) with $\rho>0$. Let $1/2 <s < (1+ \rho)/ 2$. We define the set $\vE_+$ (resp;, $\vE_-$)  of outgoing (resp., incoming) resonances of $P$ by
\[
\vE_\pm=\{\lambda \in ]0, \infty[\setminus \vT; \exists f \in \vH^{-s} \setminus \vH,  (1+ R_0( \lambda \pm i0)W)f=0\}
\]
For $\lambda \in \vE_+$, a solution to the equation $ (1+ R_0( \lambda + i0)W)f=0$ with $f \in \vH^{ -s} \setminus \vH$ is called an outgoing resonant state associated with $\lambda$. Incoming resonant states are defined similarly. $\vE_+ \cup \vE_-$ is the set of non-threshold resonances of $P$.
\end{defn}

The restriction on $s$ in the interval $]\f 1 2, \f{1+\rho} 2[$ is required by the use of Fredholm theory for compact operators in $\vH^{-s}$. The above definition seems to be $s$-dependent, but  using the properties of $R_0(\lambda\pm i0)$, one sees that if $f\in \vH^{-s}$ is a solution to the equation
$ (1+ R_0( \lambda \pm i0)W)f=0$ for some $s \in ]\f 1 2, \f{1+\rho} 2[$, then $f \in \cap_{r > \f 1 2}\vH^{-r}$. Therefore, real resonances and resonant states given in Definition \ref{defn1.1} are $s$-independent. For real resonances of non-self-adjoint Schr\"odinger operators and their role in spectral and scattering theories, the interested reader can see \cite{aw, ff, ffro, ky, sa, sch, w1, w3}. For dissipative Schr\"odinger operators, outgoing resonances are absent, but for every  positive number $\lambda>0$, one can construct a smooth compactly-supported dissipative potential such that $\lambda$ is an incoming resonance (\cite{w1}). For this reason, we do not expect $\vE_+ =\vE_-$ for the  KFP operator in general.
\\

 Denote $R(z) =( P-z)^{-1}$ for $z\not\in \sigma(P)$.  The main result of this work is the limiting absorption principle for the KFP operator with a long-range potential.\\

\begin{theorem}\label{thm1.1} Assume condition \ref{ass1} with $\rho>0$. Then

(a). The sets of accumulation points of $\sigma_p(P)$ and  $\vE_\pm$ are included in  $\vT$.

(b). Let $s >  1 /2$ . For $\lambda \in \bR_+\setminus (\sigma_p(P)\cup \vT\cup \vE_\pm)$,  the  boundary values of the resolvent
\be \label{LAP}
R(\lambda \pm i0) = \lim_{\pm \Im z >0, z\to \lambda} R(z)
\ee
exist in $\vL(s, -s)$  and are continuous in  $\lambda \in \bR_+\setminus (\sigma_p(P)\cup \vT\cup \vE_\pm)$.
\end{theorem}

We also study  decay properties of the eigenfunctions of $P$.

\begin{theorem}\label{thm1.2} Let $\lambda \in \sigma_p(P)$ and $u \in D(P)$ such that $Pu =\lambda u$.

(a). If $\lambda \in \sigma_d(P)$ and (\ref{ass1}) is satisfied with $\rho >-1$, then there exists some constant $c >0$ such that
$e^{c(\w{x} + |v|^2)}u \in \vH$.

(b). If $\lambda \in\bR_+ \setminus \vT$ and (\ref{ass1}) is satisfied with $\rho >0$, then $u$ decays polynomially: for any $s \ge 0$, $(\w{x} +\w{v})^su \in \vH$.
\end{theorem}

\begin{rmk}
Most of works in the literature on spectral analysis of the KFP operator concern the case of compact resolvent. In this connection, let us mention that  Helffer-Nier \cite{hln} conjectured that the KFP operator is of compact resolvent if and only if the associated Witten Laplacian is of compact resolvent.  See \cite{mbs, lwx, lwx2} for some partial results. Under condition (\ref{ass1}) with $\rho >-1$, the resolvent of $P$ is never compact. To prove Theorem \ref{thm1.1} and \ref{thm1.2}, we use the scattering approach initiated in \cite{w2}, considering the potential term $W$ as perturbation of the free KFP operator $P_0$.  See also \cite{nw,wz}. In \cite{w2}, the resolvent of the KFP operator is studied near the first threshold $0$ for short-range potentials ($\rho >1$) and $n=3$. It is proven that under these conditions, $0$ is not an accumulation point of eigenvalues. It remains an open question to study whether or not eigenvalues may accumulate towards a nonzero threshold for short-range potentials in dimension three.
\end{rmk}

 The remaining part of this paper  is  organized as follows. In Section 2 we recall from \cite{w2} some results on the free KFP operator and improve some of them for later use.
 In Section 3 we analyze non-threshold resonances and embedded eigenvalues of $P$ and prove that these sets do not have an accumulation point outside the thresholds $\vT$.
Once these preparations done, Theorem \ref{thm1.1} can be easily inferred (see Section 4). In Section 5 we discuss decay properties of eigenfunctions.
Appendix A. contains the proof for a technical lemma.
\\

\sect{LAP for the free KFP operator}

We first recall from \cite{w2} some results needed in this work. Let
\[
P_0 = -\Delta_v+\f{1}{4}|v|^2-\f{n}{2} + v\cdot\na_x, \quad (x,y) \in \bR^{2n},
\]
and $R_0(z) = (P_0-z)^{-1}$.
By the partial Fourier transform $\vF_{x\to \xi} $ in $x$-variables, the free KFP
$P_0$ is a direct integral of the family $\{\widehat{P}_0(\xi) ; \xi \in \bR^n\}$, where
\be
\widehat{P}_0(\xi) = -\Dl_v+\f{v^2}{4} -\f{n}{2}+ i v\cdot\xi.
\ee
The operator $\widehat{P}_0(\xi)$ is given by
\[
\widehat{P}_0(\xi) = \vF_{x\to \xi} P_0 \vF_{\xi \to x}^{-1} = -\Dl_v+\f{1}{4}\sum^n_{j=1}(v_j+2i\xi_j)^2-\f{n}{2}+|\xi|^2.
\]
 $\{\widehat{P}_0(\xi), \xi\in \bR^n\}$ is  a holomorphic family of type $(A)$ in sense of Kato with constant domain
 $D= D(-\Dl_v+\f{v^2}{4})$ in $L^2(\bR^n_v)$.  Let $F_j(s)=(-1)^je^{\f{s^2}{2}}\f{d^j}{ds^j}e^{-\f{s^2}{2}}, j \in \bN,$ be  Hermite polynomials and
$$
\varphi_j(s)=(j!\sqrt{2\pi})^{-\f{1}{2}}e^{-\f{s^2}{4}}F_j(s)
$$
the normalized Hermite functions.  For $\xi \in \bR^n$ and $\alpha=(\alpha_1, \alpha_2, \cdots, \alpha_n) \in \bN^n$, define
\be
\psi_\alpha(v) = \prod_{j=1}^n\varphi_{\alpha_j}(v_j) \mbox{ and } \psi_\alpha^\xi(v) = \psi_\alpha(v + 2i \xi).
\ee
One can check (\cite{w2}) that each eigenvalue $E_l(\xi) = l +|\xi|^2$ of $\widehat{P}_0(\xi)$  is semi-simple  and
the associated Riesz projection is given by
 \be\label{projection}
 \Pi^\xi_l\phi=\sum_{\alpha, |\alpha|=l} \langle\psi^{-\xi}_\alpha , \phi \rangle\psi^\xi_\alpha, \quad \phi \in L^2.
  \ee

\begin{prop}[\cite{w2}]\label{prop2.2}  Let $l\in \bN$ and $ l < a <l +1$ be fixed. Take  $\chi \ge 0$ and $\chi \in C_0^\infty(\bR^n_\xi)$ with supp $\chi \subset \{\xi, |\xi|^2 \le a + 4\}$, $\chi(\xi) =1$ when $|\xi|^2 \le a +3$ and $ 0 \le \chi(\xi) \le 1$.  Then one has
\be \label{R0}
\widehat R_0(z, \xi) = \sum_{k=0}^l \chi(\xi) \frac{\Pi_k^\xi}{\xi^2+k -z} + r_l(z,\xi),
\ee
for any $\xi \in \bR^n$ and $z \in \bC$ with $\Re z <a$ and $\Im z \neq 0$. Here $r_l(z, \xi)$ is holomorphic in $z$ with $\Re z < a $ verifying the estimate
\be \label{reste}
\sup_{\Re z < a, \xi \in \bR^n} \|r_l(z, \xi)\|_{\vL(L^2(\bR^n_v))} <\infty.
\ee
\end{prop}

Here and in the following, for two normed spaces $\vF$ and $\vG$,  $\vL(\vH, \vG)$ denotes the space of bounded linear operators from $\vH$ to $\vG$ and $\vL(\vH) = \vL(\vH, \vH)$. (\ref{R0}) is a kind of spectral decomposition for $\widehat R_0(z, \xi) $. The proof of the assertions on $r_l(z, \xi)$ is technical and uses time-dependent method  (cf. Proposition 2.7 and Appendix A. in \cite{w2}).\\

For a temperate symbol $a(x,v; \xi, \eta)$ (\cite{hor}), denote by $a^w(x,v, D_x, D_v)$ the associated Weyl pseudo-differential operator defined by
\bea
\lefteqn{a^w(x,v,D_x,D_v)u(x,v)} \\[2mm]
& =& \f{1}{(2\pi)^{2n}} \int\int e^{i(x-x')\cdot\xi + i(v-v')\cdot\eta} a(\f{x+x'}2, \f{v+v'}2, \xi,\eta) u(x', v') dx'dv'd\xi d\eta \nonumber
\eea
for $u \in \vS(\bR^{2n}_{x,v})$.

Remark that $ \chi(D_x) \Pi_k^{D_x}$ is a  Weyl pseudodifferential operator with a nice symbol $b_k$ independent of $x$:
\be
\chi(D_x) \Pi_k^{D_x} =  b_k^w(v, D_x, D_v)
\ee
where  $b_k(v, \xi,\eta)$ is given by
 \be
 b_k(v, \xi,\eta) = \int_{\bR^n} e^{-i v'\cdot\eta/2}\left(\sum_{|\alpha|=k} \chi(\xi) \psi_\alpha( v+ v' + 2i\xi)\psi_\alpha( v- v' + 2i\xi) \right) dv'.
 \ee
 In particular,
\be
b_0(v, \xi,\eta) = 2^{\f n2}\chi(\xi)e^{-v^2 -\eta^2 + 2iv\cdot\xi + 2\xi^2}.
\ee
These pseudodifferential operators  belong to $\vL(r,s; r',s)$ for any $r, r', s \in \bR$  (Ch. 18, \cite{hor}). Set
\[
 r_l(z) = \vF_{\xi\to x}^{-1} r_l(z, \xi)\vF_{x\to\xi}
 \]
 Proposition \ref{prop2.2} implies the following representation formula for the free resolvent
\be \label{RR0}
 R_0(z) = \sum_{k=0}^l b_k^w(v, D_x, D_v)(-\Delta_x +k -z)^{-1} + r_l(z),
\ee
for  $z \in \bC$ with $\Re z <a$ and $\Im z \neq 0$ and
 that  $ r_l(z)$ is  bounded on $L^2$ and  holomorphic in $z$ with $\Re z <a$, $a \in ]l, l+1[$.
 \\

 \begin{lemma} \label{lem2.3} Let $z\in \bC\setminus \bR_+$. For any $m\in \bN$, $ s \in \bR$ and $r\in[0, 2m]$, one has $R_0(z)^m \in \vL(-r, s; 2m-r,s)$.
 \end{lemma}
\pf Let $m=1$. Consider first the case $r=0$ and $s\in \bN$. For $s=1$, one has
\[
\w{x}R_0(z)\w{x}^{-1} = R_0(z)\left( 1+ [-v\cdot\nabla_x, \w{x}] R_0(z) \w{x}^{-1}\right).
\]
The subelliptic estimate (\ref{subelliptic}) shows that $R_0(z) \in\vL(0,2)$ and  $[-v\cdot\nabla_x, \w{x}]R_0(z) \in \vL(0,0)$. Therefore, $\w{x}R_0(z)\w{x}^{-1} \in\vL(0, 2)$.
Similarly, using the formula
 \[
 \w{x}^k R_0(z)\w{x}^{-k}  = \w{x}^{k-1} R_0(z)\left( \w{x}^{-k+1} + [-v\cdot\nabla_x, \w{x}] R_0(z) \w{x}^{-k}\right),
 \]
we can show by an induction that $\w{x}^k R_0(z)\w{x}^{-k} \in\vL(0, 2)$ for all $k \in \bN$. This proves  $R_0(z) \in \vL(0,k; 2, k)$.
By an argument of duality, one has $R_0(z)^*\in \vL(-2, -k; 0, -k)$ for $k \in \bN$.
Note that $P_0^* = J P_0 J$ where $ J f(x, v) = f(x,-v)$ and $J$ conserves any $\vH^{r,s}$ norm. From $R_0(\overline{z})^* = J R_0(z)J$,  we deduce that $R_0(z)\in\vL(-2,-k; 0;-k)$ for $k \in \bN$. A complex
interpolation between $ (r_1, s_1) =(-2, -k)$ and $(r_2, s_2) =(2, k)$ gives that $R_0(z) \in \vL(-r, s; 2-r, s)$ for $r\in[0, 2]$ and  $ |s| \le k$, for any $k \in \bN$.
The lemma is proven for $m=1$.
 \\

  When $m \ge 2$, one can use the result obtained for $R_0(z)$  and the above commutator method to show that $R_0(z)^m  \in \vL(0,s; 2m, s)$ for any $s\in \bR$. The details are omitted here. The arguments of duality and complex interpolation used above allow us to conclude that   $R_0(z)^m \in \vL(-r, s; 2m-r,s)$ for $s\in \bR$ and $r\in[0, 2m]$.
\ef

\begin{prop}\label{prop2.1} Let $n \ge 1$. For any $s>\f 1 2$,
the boundary values of the resolvent
\[
R_0(\lambda \pm i\epsilon) =  \lim_{ \pm \Im z >0, z \to \lambda}R_0(z)
\]
 exist
in $\vL(0,s; 2,-s)$ for $\lambda \in \bR_+\setminus \vT$ and is continuous in $\lambda$. More generally, let $k \in \bN$ and $s> k + \f 1 2$ and $ r, r' \ge 0 $ with $r+r' \le 2(k+1)$, $R_0(\lambda \pm i0)$ is of class $C^k$ in $\vL(-r,s; r',-s)$ for $\lambda \not\in \vT$.
\end{prop}
\pf When $k=0$, $r=0$ and $r'=2$, the result has been proven in \cite{w2} by using  (\ref{RR0}) and the known results for $(-\Delta_x -(\mu \pm i0))^{-1}$ in weighted $L^2$-spaces(cf. \cite{ag1}).
In the general case, the existence and the continuity of $\f{d^k}{d\lambda^k}R_0(\lambda \pm i0)$ in $\vL(0, s; 0, -s)$, $s> k +\f 1 2$, can be deduced by the same argument.
To see its continuity in weighted Sobolev spaces, we use Lemma \ref{lem2.3} to conclude $R_0( -1)^{k+1}  \in \vL(0, s; 2k+2, -s)$ and the equation
\[
R_0(\lambda \pm i0)^{k+1} = R_0( -1)^{k+1} \left( 1  + (\lambda+1) R_0(\lambda \pm i0)\right)^{k+1}
\]
 to show that  $\f{d^k}{d\lambda^k}R_0(\lambda \pm i0)$ is continuous in $\lambda$ in topology of $\vL(0, s;  2k +2, -s)$, if $s> k +\f 1 2$. The general statement of Proposition \ref{prop2.1} follows from the arguments of duality and complex interpolation already used in Lemma \ref{lem2.3}.
\ef

\sect{Real resonances and embedded eigenvalues}

Let us now consider the full KFP operator $P$ with a long-range potential satisfying the condition \ref{ass1} with $\rho>0$.

\begin{lemma}\label{lem3.1}
  Let $\lambda_0 \in \bR_+ \setminus \vT$ be an embedded eigenvalue of $P$ and $u$ an associated eigenfunction: $P u=  \lambda_0 u$. Then one has
\be \label{e3.1}
u= - R_0(\lambda_0 \pm i0)Wu.
\ee
\end{lemma}
\pf Proposition \ref{prop2.1} shows that $R_0(\lambda_0 \pm i0) Wu \in \vH^{1, -s}$, $1/2 <s < (1+\rho)/2$. To show (\ref{e3.1}) for sign $+$, we set $z= \lambda_0 + i\epsilon $ with $\epsilon >0$ and write $u$ as
\[
u = -R_0(z) Wu - i \epsilon R_0(z)u.
\]
To see that $\lim_{\epsilon \to 0_+} \epsilon R_0(z)u =0$ in $\vH$, we use (\ref{RR0}) with $l > \lambda_0$. Clearly
\[
\lim_{\epsilon \to 0_+} \epsilon r_l(z)u =0,
\]
 because $r_l(w)$ is holomorphic in $w$ with $\Re w < l +a$ for some $a>0$.  The symbol of $\epsilon (-\Delta + k -z)^{-1}$ is equal to
 $\epsilon (|\xi|^2 + k - \lambda_0 - i \epsilon)^{-1}$ which converges pointwise to $0$ and  is uniformly bounded by $1$. The Dominated Convergence Theorem shows that
\[
\lim_{\epsilon \to 0_+}  \sum_{k=0}^l \epsilon b_k^w(v, D_x, D_v)(-\Delta_x +k -z)^{-1}u = 0
\]
 in $\vH$.  This proves that $\lim_{\epsilon \to 0_+} \epsilon R_0(z)u =0$ in $\vH$, hence $u = - R_0(\lambda_0 + i0) Wu$ in $\vH$. Similarly one can show that $u= - R_0(\lambda_0-i0)Wu$.
\ef

\begin{lemma}\label{lem3.2} Let $\lambda_0\in \bR_+\setminus \vT$ and $ 1/2 <  s < {1+\rho}/2$. If $u \in \vH^{-s}\setminus\{0\}$ satisfies both equations
\[
u = - R_0(\lambda_0 \pm i0) Wu,
\]
then $u \in \vH$ and $\lambda_0 \in \sigma_p(P)$.
 \end{lemma}
\pf Making use of (\ref{RR0}) with $l> \lambda_0$, one has
\be \label{u}
u =  -\sum_{k=0}^l b_k^w(v, D_x, D_v)(-\Delta_x +k -(\lambda_0 \pm i0))^{-1}Wu - r_l(\lambda_0)Wu.
\ee
Since $Wu \in \vH^{-1, 1+\rho -s}$ and $1+\rho -s > s$, $ r_l(\lambda_0)Wu\in \vH$.   Let $d =d(\lambda_0, \vT)$. Take $\chi_1 \in C_0^\infty(\bR^n_\xi \setminus 0)$ such that
\[
\chi_1(\xi) =1 \quad \mbox{ for } \f d 2 \le |\xi|^2 \le l +1.
\]
Let $\chi_2(\xi) = 1-\chi_1(\xi)$ and decompose
\[
u = u_1 + u_2 \mbox{ with } u_j = \chi_j(D_x)u.
\]
 From (\ref{u}) and the fact that  $\chi_2(\xi) (|\xi|^2+ k -(\lambda_0 \pm i0))^{-1}$, $k=0, \cdots, l$, are bounded for $\xi \in \bR^n$, we deduce that $u_2 \in \vH$. To show that $u_1 \in \vH$, we decompose $\chi_1(D_x)$ as
\[
\chi_1(D_x) = B_+ + B_-
\]
where $B_\pm = b_\pm^w(x, D_x)$ are Weyl pseudo-differential operators with bounded symbols $b_\pm(x, \xi)$ such that
\[
\mbox{ \rm supp }  b_\pm \subset \{(x, \xi); |x|\le R\} \cup \{(x, \xi); \pm \widehat x \cdot \widehat \xi > (1-\delta_\pm)\}
\]
for some $R>1$ and $0<\delta_\pm <1$. Here $\widehat x = \f x{|x|}$, $\widehat \xi =\f \xi {\xi}$ for $x, \xi \neq 0$. One has the microlocal resolvent estimates for $-\Delta_x$ in $L^2(\bR^n_x)$ (\cite{ik}): for $\mu >0$ and $s >  1 /2$,
\be \label{microl}
\| \w{x}^{s-1}B_\mp (-\Delta_x - (\mu \pm i0))^{-1}\w{x}^{-s} \| <\infty.
\ee
It follows that
\be \label{Bpm}
B_\mp b_k^w(v, D_x, D_v)(-\Delta_x +k -(\lambda_0 \pm i0))^{-1} \in \vL(r, s; r', 1-s)
\ee
for all $r,r' \in \bR$ and $s>1/2$.
Write $u_1$ as
\bea
\lefteqn{ u_1 = B_- u+ B_+ u} \nonumber \\
 & = & - B_- \left(\sum_{k=0}^l b_k^w(v, D_x, D_v)(-\Delta_x +k -(\lambda_0 + i0))^{-1}Wu + r_l(\lambda_0)Wu\right) \label{u1} \\
 &  & -  B_+\left(\sum_{k=0}^l b_k^w(v, D_x, D_v)(-\Delta_x +k -(\lambda_0 -i0))^{-1}Wu + r_l(\lambda_0)Wu\right). \nonumber
\eea
Applying (\ref{microl}) with $r= 1+ \rho -s$ and $\mu =\lambda -k$, $k=0, \cdots, l$, one obtains that $u_1 \in \vH^{-(s-\rho)}$. If $s-\rho \le 0$, then  $u = u_1 + u_2 \in \vH$.
The proof is finished.
If $s-\rho >0$, then $u \in \vH^{-(s-\rho)}$.  We apply (\ref{microl}) once more with $r = 1 +\rho -(s-\rho)$ and deduce from (\ref{u1}) that $u_1 \in \vH^{-(s-2\rho)}$ which implies $u \in \vH^{-\max\{0, s-2\rho\}}$.  Repeating these arguments for at most a finite number of times, we obtain $u \in \vH$.
\ef

\begin{rmk} We can also give large-$x$ expansions of resonant states. If $u$ is an outgoing resonant state  associated with a resonance $\lambda\not\in \vT$, then the argument used above shows that
\[
u- \sum_{k=0}^{l_0} b_k^w(v, D_x, D_v)(-\Delta_x +k -(\lambda + i0))^{-1}Wu \in \vH,
\]
where $l_0$ is the integer such that $l_0 < \lambda < l_0 +1$. We can study the large $|x|$ behavior of $u$, by making use of the  distributional kernel of $(-\Delta_x -z)^{-1}$ given by
\[
K(x, y; z) = \f{i} 4 \left( \f{\sqrt z}{2|x-y|}\right)^{\f n 2 -1} H^{(1)}_{n/2-1}(\sqrt z |x-y|), \quad x, y \in \bR^n,
\]
where $\Im \sqrt z \ge 0$ and $H^{(1)}_\nu$ is the Bessel function of the third kind (\cite{ol}).  In particular, for $n=3$, one has simply
\[
K(x, y; z) = \f{e^{i\sqrt z |x-y|}}{4\pi |x-y|}.
\]
For $\rho>1$ and $n=3$, one can check that
\[
(-\Delta_x +k -(\lambda + i0))^{-1}Wu  - \phi_k \in \vH,
\]
where
\be
\psi_k(x,v) = \f{e^{i\sqrt {\lambda-k} |x|}}{4\pi |x|} \int_{\bR^3_y} (Wu)(y,v) dy.
\ee
This proves
\be
u = \sum_{k=0}^{l_0} b_k^w(v, D_x, D_v) \psi_k  + h,
\ee
with $h \in \vH$. In particular, if
\be
\int_{\bR^3} \nabla V(y) u(y,v) dy =0, \quad  \mbox{ a.e. in } v,
\ee
then $u  $ is an eigenfunction and $\lambda$  an embedded eigenvalue of $P$.  Similar results are also true for incoming resonant states.
\end{rmk}

Let $K \subset \overline{\bC}_\pm$ and  $1/ 2 < s <(1+\rho)/ 2$. Let $E_K$ be the subspace of $ \vH^{-s}$ spanned by eigenfunctions and outgoing or incoming resonant
states (according to sign $+$ or $-$) associated with eigenvalues and resonances in $K$.

\begin{prop} \label{prop3.3} If $K \subset \overline{\bC}_\pm$ is compact and $K\cap \vT =\emptyset$, then $ E_K$ is of finite dimension. In particular, the sets of accumulation points of $\sigma_p(P)$ and $\vE_\pm$ are included in $ \vT$.
\end{prop}
\pf To fix the idea, let $K \subset \overline{\bC}_+$ and $ 1/ 2 < s' < s$.  Since $ 1+\rho - s >  1/ 2$, Proposition \ref{prop2.1}
shows that
\[
\| R_0(z)W\|_{\vL(0, -s; 1, -s') } \le C_K
\]
uniformly for $z\in K$. Here $R(z)$ is understood as $R(\lambda + i0)$ if $z =\lambda \in K\cap\bR$. Let $B$ be the unit ball in $E_K$ and $u \in B$. Using the equation
$u= - R_0(z_0)Wu$ for some $z_0\in K$, one has
\[
\|u\|_{\vH^{1, -s'} }\le C_K \|u\|_{\vH^{-s}} = C_K,
\]
 This means that $B$ is bounded in $\vH^{1, -s'}$. Since the canonical injection $\vH^{1, -s'} \hookrightarrow \vH^{-s}$ is compact for $1/2 < s'< s$, $B$ is compact in $\vH^{-s}$.
 This proves that $\dim E_K <\infty$ by F. Riesz's criteria on the finiteness of dimension for normed vector spaces.
\ef

Proposition \ref{prop3.3} implies that the geometric multiplicities of embedded eigenvalues of $P$ are finite, but we have no information on their algebraic multiplicities.
\medskip

\sect{LAP for the perturbed KFP operator}

{\bf \noindent Proof of Theorem \ref{thm1.1}.} Part (a) of Theorem \ref{thm1.1} is a consequence of Proposition \ref{prop3.3}.
 For Part (b), we only study $R(\lambda + i0)$.  Let $  1/2 < s < (1+\rho)/ 2$. For $\lambda \in \bR_+ \setminus (\vT \cup \sigma_p(P) \cup \vE_+)$, $ 1 + R_0(\lambda+i0)W$ is injective in $\vH^{-s}$ and $R_0(\lambda)W$ is a compact operator in $\vH^{-s}$. Therefore $1+ R_0(\lambda + i0)W$ is invertible in $\vL(-s, -s)$.  By the continuity of $R_0(\lambda + i0)$ (Proposition \ref{prop2.1}),
$(1 + R_0(\mu + i0) W)^{-1}$ also exists in $\vL(-s, -s)$ for $\mu$ near $\lambda$ and is continuous in $\mu$. Taking the limit $z \to \mu + i0$ in the equation
\[
R(z) = (1 + R_0(z)W)^{-1} R_0(z)
\]
we see that for $\mu$ near $\lambda$, $R(\mu + i0) = \lim_{z \to \mu + i0} R(z)$ exists in $\vL(s,-s)$ and is given by
\be
R(\mu + i0) = (1 + R_0(\mu + i0)W)^{-1} R_0(\mu + i0).
\ee
As  $\vL(s,-s)$-valued function,  $R(\mu + i0)$ is continuous in $\mu$, because $R_0(\mu + i0) \in \vL(s,-s)$ is continuous in $\mu$ and the same is true for $((1 + R_0(\mu + i0)W)^{-1}$ as $ \vL(-s, -s)$-valued function. This proves Part (b) for $R(\lambda +i0)$. The proof for $R(\lambda -i0)$ is the same.
\ef

\begin{rmk}
If (\ref{ass1}) is satisfied with $\rho> 2k$, $k \in\bN^*$, then one can use Proposition \ref{prop2.1} to show that $R(\lambda\pm i0)$ is of class $C^k$ in $\vL(0, s; 2k+2, -s)$, $s> k + \f 1 2$, for $\lambda \in \bR\setminus (\vT \cup \sigma_p(P) \cup \vE_\pm)$.
\end{rmk}

\sect{Decay properties of eigenfunctions}

We prove Theorem \ref{thm1.2} by studying separately discrete and embedded eigenvalues.   The following result is a cosmetic improvement of Theorem \ref{thm1.2} (a).
\\

\begin{theorem} \label{thm5.1} Suppose that condition (\ref{ass1}) is satisfied for some $\rho >-1$. Then there exist some positive constant $c_0$ and positive-valued function $\tau(\zeta)$ (defined by (\ref{tau}) below) for $\zeta\in \bC\setminus \bR$ such that if $u$ is an eigenfunction of $P$ associated with an eigenvalue $z\in\sigma_d(P)$, then $e^{c_0 \tau(z) (\w{x} + |v|^2)} u \in \vH$.
\end{theorem}
\pf Let $u \in D(P)$ and $z\in \sigma_d(P)$ such that that $P u = z u$.
Then $z\not \in \bR_+$ and one has
\be \label{e5.1}
u = -R_0(z)Wu.
\ee

Let $\chi (s)$ be a smooth real-valued function on $\bR$ such that
\[
\chi(s) = \left\{ \begin{array}{ccc}
s, & \quad &  \mbox{ if } \quad s\in [0, 1], \\
\f 3 2,  & \quad &  \mbox{ if } \quad s\ge 2.
\end{array} \right.
\]
Let $\chi_r(s) = r \chi (\f s r)$ for $r \ge 1$. Then $\chi_r(s) \to s$ pointwise  as $r\to \infty$ and $|\chi'_r(s)| \le C $, uniformly for $r \ge 1$ and $s\ge 0$.
Denote
\[
\phi_r(x, v) = \chi_r(\w{x}) + \chi_r(|v|^2).
\]
Set $P_{0, r} = e^{a \phi_r} P_0 e^{-a \phi_r}$ with $a>0$ which is to be chosen later.  One can compute
\[
P_{0,r} = P_0 + Q_r,
\]
where
\[
Q_r =  a \left(4 \chi_r'(|v|^2) v\cdot \nabla_v + \Delta_v(\chi_r(|v|^2)) - a|\nabla_v \chi_r(|v|^2)|^2 - \chi_r'(\w{x}) v\cdot \f{x}{\w{x}}\right)
\]
By the choice of $\chi_r$, there exists some constant $C$ such that
\[
\|Q_r R_0(z)\| \le C a (\|  v\cdot \nabla_v R_0(z)\| + \| v\cdot \f{x}{\w{x}} R_0(z)\|) + (1+a)\|R_0(z)\|).
\]
Define $\tau(\zeta)$ by
\be \label{tau}
\tau(\zeta) =\min\{ (\|  v\cdot \nabla_v R_0(\zeta)\| + \| v\cdot \f{x}{\w{x}} R_0(\zeta)\| +\|R_0(z)|)^{-1}, \|R_0(\zeta)\|^{-\f 1 2}\}
\ee
for $\zeta \in \bC\setminus \bR_+$.
For $ a=  c_0 \tau (z)$,  one has
\[
\|Q_rR_0(z)\| \le C (c_0 + c_0^2) <1
\]
if $c_0 >0$ is chosen appropriately small. With such a choice of $a$,  the resolvent $ R_{0, r}(z) =(P_{0,r}-z)^{-1}$ exists  for all $r\ge 1$ and is given by
\[
R_{0, r}(z) = R_0(z) (1 + Q_rR_0(z))^{-1}.
\]
It follows that there exists some constant $C$ depending on $z$ such that
\be
\|R_{0, r}(z)\| \le C
\ee
uniformly for $r \ge 1$. From the subellipticity of $P_0$, one deduces that $R_{0,r}(z) \in \vL(-s,0; 2-s,0)$, $s\in[0, 2]$, is uniformly bounded for $r\ge 1$.
In particular, $\|R_{0,r}(z)(\w{v} + \w{D_v})\| \le C $ uniformly for $r \ge 1$.
 \\

 Let $u_r = e^{c_0 \tau(z) \phi_r}u$. Equation (\ref{e5.1}) gives
\be
u_r = - R_{0,r}(z)(\nabla_v  \chi_r(|v|^2)  -\nabla_v)\cdot \nabla_xV u_r
\ee
Since $\nabla_x V = O(\w{x}^{-1-\rho})$ and $1+ \rho >0$,  there exists some constant $C$ independent of $R_0$ and $r$ such that for any $R_0>1$,
\[
\|\nabla_x V u_r\| \le C_{R_0} \|u\| + C \w{R_0}^{-1-\rho}\|u_r\|.
\]
Seeing that there exists some constant $M>0$ such  that
\[
 \| R_{0, r}(z) (\nabla_v  \chi_r(|v|^2)  -\nabla_v)\| \le C' \|R_{0, r}(z)\|_{\vL(-1,0; 0,0)} \le M,
 \]
 uniformly for $r \ge 1$,   one has
\beas
\|u_r\| &\le &    \| R_{0, r}(z) (\nabla_v  \chi_r(|v|^2)  -\nabla_v)\| \; \|\nabla_x V u_r\|  \\
& \le & M ( C_{R_0} \|u\| + C \w{R_0}^{-1-\rho} \|u_r\|), \quad \forall r \ge 1.
\eeas
Taking $R_0$ appropriately large so that  $M C \w{R_0}^{-1-\rho} < \f 1 2$, one obtains
\[
\| u_r\|  \le 2 M C_{R_0} \|u\|,  \quad r \ge 1.
\]
Since $u_r(x, v) \to e^{c_0\tau(z)(\w{x} + v^2)} u(x, v)$ pointwise as $r\to \infty$, Fatou's Lemma shows that $e^{c_0\tau(z)(\w{x} + v^2)} u \in \vH$.
\ef

\begin{rmk} Theorem \ref{thm5.1} gives some quantitative information about  decay properties of eigenfunctions.  Comparing with the Green function of $P_0$, we believe that the phase function $\w{x} + |v|^2$ is the appropriate one to use in this context, but the numerical factor $c_0 \tau (z)$ before it can very likely be improved.
 If $z = b + i \mu$ with $b>0$ fixed and $\mu \in \bR$, one has uppers bounds for
$
\|  v\cdot \nabla_v R_0(z)\| $, $ \| v\cdot \f{x}{\w{x}} R_0(z)\|)$ and  $\|R_0(z)\| $
for $|\mu| \neq 0$. However to calculate $\tau(z)$, one also needs to evaluate their lower bounds, which seems not to be an easy task.
\end{rmk}

Theorem \ref{thm1.2} (b) is proven by applying the boost-up argument of Lemma \ref{lem3.2} to
\[
u = -R_0(\lambda \pm i0) Wu
\]
 to improve the decay of $u$. To do this, we need to show that the remainder $r_l(\lambda)$ given by (\ref{RR0}) with $l > \lambda$ and $z=\lambda$ conserves  polynomial decay. Since $\nabla_xV(x)$ does not operate in $\vF^{r,s}$ for $r$ large, we modify the definition of $\vF^{r,s}$ and introduce $\vG^{r,s}$ by
\[
\vG^{r,s} =\{f \in \vS'(\bR^{2n}; \w{ 1 -\Delta_v + v^{2}}^{\f r 2} \w{x}^s f \in L^2\}.
\]
$\vG^{r, s}$ is equipped with its natural norm. Then, $\nabla_xV \in \vL( \vG^{r,s}; \vG^{r, s+1 +\rho})$ for all  $r,s$. The proof of the following technical lemma is given in Appendix A.

 \begin{lemma} \label{lem5.3}  For any $r, s \in \bR$, one has $r_l(\lambda)\in \vL(\vG^{r, s}; \vG^{2+r,s})$.
 \end{lemma}

 \medskip

{\bf \noindent Proof of Theorem \ref{thm1.2} (b).} As in the proof of Lemma \ref{lem3.2}, we  apply (\ref{RR0}) with $l> \lambda$  to $u = -R_0(\lambda \pm i0) Wu$ and write
\[
u =  -\sum_{k=0}^l b_k^w(v, D_x, D_v)(-\Delta_x +k -(\lambda \pm i0))^{-1}Wu - r_l(\lambda)Wu.
\]
Since $u \in \vH$ and $Wu \in \vG^{-1, 1 +\rho}$,  Lemma \ref{lem5.4} shows that
\[
r_l(\lambda)Wu \in \vG^{1, 1 +\rho}.
\]
Introduce the microlocal partition of unity $\chi_2(D_x) + B_+ + B_- =1$ as constructed in Lemma  \ref{lem3.2}.
Then an argument of support in $\xi$ gives
\[
\chi_2(D_x) b_k^w(v, D_x, D_v)(-\Delta_x +k -(\lambda \pm i0))^{-1}Wu \in \vG^{r, 1+\rho}
\]
and (\ref{microl}) with $r =1+\rho_0$  shows that
\[
B_\mp b_k^w(v, D_x, D_v)(-\Delta_x +k -(\lambda \pm i0))^{-1}Wu \in \vG^{r, \rho}
\]
for all $r \in \bR$. Putting these pieces together, one obtains $u \in \vG^{1, \rho}$.
Now $Wu \in \vG^{0, 1+ 2\rho}$. Repeating the above arguments with this improved property on $Wu$, we obtain $u \in \vG^{2, 2\rho}$.
By an induction on $m$, we can show that $u \in \vG^{m, m\rho}$ for any $m\in \bN$. Consequently, $u$ decays polynomially.
\ef

\begin{rmk} For  $N$-body self-adjoint Schr\"odinger operator $H$ with long-range interactions, it is known (\cite{ag2,fh}) that  an eigenfunction  associated with a non-threshold eigenvalue $\lambda$ of $H$ decays like $O(e^{-(1-\epsilon) \sqrt{d(\lambda)}\: |x|})$ for any $\epsilon >0$ and small, where $d(\lambda)$ is the distance between $\lambda$ and the thresholds of $H$ located on the right-hand side of $\lambda$. For the KFP operator $P$ with a long-range potential, one may naturally ask the questions whether there exists  an exponential decay of eigenfunctions associated with an embedded eigenvalue $\lambda \in \bR_+\setminus \vT$  and whether their decay rate can be described in terms of the distance between $\lambda$ and $\vT$.
\end{rmk}

\appendix

\sect{Proof of Lemma \ref{lem5.3}}

We first give two results needed for the proof of Lemma \ref{lem5.3}. The first one concerns the evolution of elementary observables by the semigroup $e^{-t P_0}$.

\begin{lemma}[Lemma 2.12, \cite{w2}] \label{lemA.2} Let $n \ge 1$. For $t\ge 0$ and $0 \le s \le t$, one has the following equalities as operators from $\vS(\bR^{2n}_{x,v})$ to $L^2(\bR_{x,v}^{2n})$:
\bea
e^{-(t-s)P_0}v_j e^{-sP_0} &=& e^{-t P_0} (v_j \cosh s - 2  \partial_{v_j} \sinh s + 2 (\cosh s -1) \partial_{x_j}) \label{v1} \\
e^{-(t-s)P_0}\partial_{v_j }e^{-sP_0} &=& -\f 1 2 e^{-t P_0} (  (v_j \sinh s - 2  \partial_{v_j} \cosh s   + 2 \partial_{x_j} \sinh s
 ))  \label{v2} \\
e^{-(t-s)P_0}x_j e^{-sP_0} &=& e^{-t P_0} (x_j+ v_j \sinh s - 2 (\cosh s -1)\partial_{v_j} \nonumber \\
& &  + 2 (\sinh s -s) \partial_{x_j})  \label{x1}
\eea
\end{lemma}

The second one is about the smoothing properties of some integrals of $e^{-tP_0}$.

\begin{lemma}\label{lem5.4}  Let $k \in \bN$. If $\theta(t)$ is a smooth function with  $\theta^{(j)}(0) =0$
for $j=0, \cdots, k$, then
for any fixed $T>3$, the operator
\be \label{Ik}
(P_0+ 1)^{k+1}\int_0^T e^{-tP_0} \theta(t)\; dt \in \vL(\vH).
\ee
\end{lemma}
\pf  For an abstract $C_0$-contraction semigroup generated by a maximally accretive operator, results similar to (\ref{Ik}) hold true for $k =0$, but in general they fail for $k \ge 1$.  To prove (\ref{Ik})  for all $k \in \bN$,  we use  the following formula  (cf. (A.10) in \cite{w2}):
 \be \label{eqA.2}
\sum_{l=0}^{\infty}e^{-t(l+|\xi|^2)}\|\Pi_l^\xi\|
 \le \f{e^{-|\xi|^2(t-2 - \f{4}{e^t-1}) }}{(1-e^{-t})^n},
 \ee
 for $t>0$ and $\xi\in \bR$. Here $\Pi_l^\xi$ is the Riesz projection of $\widehat P_0(\xi)$ with eigenvalue $l +|\xi|^2$.  For any $k\in\bN$, one has
 \beas
\lefteqn{ \| (\widehat P_0(\xi) + 1)^k e^{-t \widehat P_0(\xi)}\|}\\
 & = & \| \sum_{l=0}^{\infty}(l+|\xi|^2)^k e^{-t(l+\xi^2)}\Pi_l^\xi\| \\
&\le &    \sum_{l=0}^{\infty}\|(l+|\xi|^2)^k e^{-t(l+\xi^2)}\Pi_l^\xi\| \\
& \le  & C_{k, \delta} \sum_{l=0}^{\i}e^{-(t-\delta)(l+|\xi|^2)}\|\Pi_l^\xi\|  \quad \mbox{ ($ 0<\delta <t$)}\\
&\le &  C_k e^{-(t-T)|\xi|^2},
\eeas
for $\xi \in \bR^n$ and $t \ge T>3$ (with $\delta >0$ chosen appropriately). Consequently  for any $k \in \bN$, one has
\be \label{uniformk}
\| (P_0 + 1)^k e^{-tP_0}\| \le C_k,
\ee
 uniformly  for $t \ge T >3$.
Lemma  \ref{lem5.4} follows from integration by parts and the above estimate with $t =T>3$.
\ef

 Lemma \ref{lem5.4}  shows  that if  $\theta^{(j)}(0) =0$
for $j=0, \cdots, k$, then
\be \label{eA.1}
\int_0^T e^{-tP_0} \theta(t) dt \in \vL(\vF^{r,0}, \vF^{r + 2(k+1), 0})
\ee
 for $r \in[-2(k+1), 0]$. In particular, $\int_0^T e^{-tP_0} \w{D_x}^{\f{2(k+1)}{3}} \theta(t) dt\in \vL(\vH)$. \\

{\bf \noindent Proof of Lemma \ref{lem5.3}.} Let $\chi_1= 1-\chi$, where $\chi$ is the cut-off used in Proposition \ref{prop2.2}. By the partial Fourier transform in $x$, $r_l(\lambda)$ is the multiplication by $r_l(\lambda, \xi)$ given by (2.40) in \cite{w2}. We decompose it as
 \be
 r_{l}(\lambda, \xi) =  S_1(\lambda, \xi) +  S_2(\lambda, \xi) + S_3(\lambda, \xi),
 \ee
 where for some constant $T >3$,
 \beas
 S_1(\lambda, \xi) & = & \int_T^\infty \left( \chi_1(\xi) e^{-t(\widehat P_0(\xi)-\lambda)} + \chi(\xi) (\sum_{k=l+1}^\infty e^{-(k +|\xi|^2)t}\Pi_k^\xi) \right)\; dt,\\
 S_2(\lambda, \xi)  &= &   \int_0^T  \chi_1(\xi) e^{-t(\widehat P_0(\xi)-\lambda)} \; dt,  \\
 S_3(\lambda, \xi)  &= & \int_0^T   \chi(\xi) e^{t \lambda}\left(e^{-t \widehat P_0(\xi)} -\sum_{k= 0}^l e^{-(k +|\xi|^2)t}\Pi_k^\xi \right ) \; dt.
 \eeas
 Denote
 \[
 S_j =S_j(\lambda, D_x) = \vF_{\xi\to x }^{-1}S_{j}(\lambda, \xi) \vF_{x\to\xi},
 \]
 for $j=1,2, 3$. \\

 $S_1(\lambda, \xi)$ is smooth and rapidly decreasing in $\xi$. Therefore $S_1\in \vL(\vG^{0, s}; \vG^{0,s})$ for any $s\in \bR$.
 From (\ref{uniformk}) and the subelliptic estimate for $P_0$, it follows that $ S_1 \in  \vL(\vG^{r, s}; \vG^{r',s})$ for any $r, r', s\in \bR$.\\

 To study $S_2$, we use the method of commutators. (\ref{x1}) with $s=0$ gives
 \be
 [x_j, e^{-tP_0} ]= e^{-t P_0} ( v_j \sinh t - 2 (\cosh t -1)\partial_{v_j}  + 2 (\sinh t -t) \partial_{x_j})
 \ee
for $j=1, \cdots, n$. Hence
\bea
\lefteqn{[S_2, x_j]}  \nonumber \\
 &= & -i(\partial_{\xi_j}\chi_1)(D_x) \int_{0}^T  e^{-t(P_0 -\lambda)} dt  \label{rl1} \\
& & -
\chi_1(D_x) \int_{0}^T e^{-t(P_0-\lambda)} (v_j \sinh t - 2 (\cosh t-1)\partial_{v_j}
 + 2 (\sinh t - t) \partial_{x_j}) dt  \nonumber
\eea
Seeing that
\[
\sinh t = O(t), \quad  \cosh t-1 = O(t^2), \quad  \sinh t - t = O(t^3),
\]
as $t \to 0$, we can apply Lemma \ref{lem5.4} to the three terms in the last integral in (\ref{rl1}) with, respectively, $k=0$, $1$ and $2$. It follows that
 \[
  [S_2, x_j] \in \vL(\vG^{0, 0}; \vG^{2,0}),
 \]
 for $j=1, \cdots, n$.  This proves
 \[
 x_j S_2\w{x}^{-1} =  ( S_2 x_j - [S_2, x_j] )\w{x}^{-1}  \in \vL(\vG^{0, 0}; \vG^{2,0}).
 \]
 Consequently,  $S_2 \in \vL(\vG^{0, 1}; \vG^{2,1})$.
  By an induction on $k$ and using the commutator method,  one can show
 \[
 S_2\in \vL(\vG^{0, k}; \vG^{2,k}), \quad \forall  k\in \bN.
 \]
 By the argument of duality and complex interpolation, we deduce that
\[
  S_2\in
 \vL(\vG^{0, s}; \vG^{2,s}), \quad s\in \bR.
 \]

 Similarly, using (\ref{v1}) and (\ref{v2}), we can evaluate the commutators
 $[S_2, v_j]$ and $[S_2, \partial_{v_j}]$.  For $[S_2, v_j] $,  one has
 \bea
\lefteqn{[S_2, v_j]}  \label{S2vj} \\
 &= &
\chi_1(D_x) \int_{0}^T e^{-t(P_0-\lambda)} (v_j (\cosh t -1)  - 2 \partial_{v_j} \sinh t  + 2 (\cosh t -1) \partial_{x_j}) dt
   \nonumber
\eea
 Applying Lemma \ref{lem5.4} and noticing that $(P_0+1)^{-1}v_j$, $(P_0+1)^{-1}\partial_{v_j}$ and $(P_0+1)^{-2}\partial_{x_j}$ belong to $\vL(s,s)$ for any $s$ (cf. Lemma \ref{lem2.3}), we infer that $[S_2, v_j] \in \vL(\vG^{0,s}, \vG^{2,s})$. Similarly, we can prove $[S_2, \partial_{v_j}] \in \vL(\vG^{0,s}, \vG^{2,s})$, using (\ref{v2}). We deduce that
 $
 S_2 \in \vL(\vG^{1,s}, \vG^{3,s}).
 $
 By an induction, we can show that $S_2 \in \vL(\vG^{k,s}, \vG^{k+2,s})$ and conclude that $S_2 \in \vL(\vG^{r,s}, \vG^{r+2,s})$ for $r,s \in \bR$. \\

The term $S_3$ is easier to study. Since $\chi(\xi)$ is of compact support  in $\xi$, $S_3$ is a convolution with a rapidly decaying function in $x$ variable, which gives that $S_3 \in \vL(\vG^{0,s}, \vG^{2,s})$  for $s\in \bR$. In addition, for any $k \in \bZ$, $ \chi(D_x) (P_0 +1)^k ( 1 - \Delta_v + |v|^2)^{-k} $ is bounded.  We conclude from Lemma \ref{lem2.3} and the equation
\beas
\lefteqn{( 1 - \Delta_v + |v|^2)^{-k} S_3 ( 1 - \Delta_v + |v|^2)^{k} }\\
& =& \left(( 1 - \Delta_v + |v|^2)^{-k} (P_0 +1)^k\right) S_3 ((P_0 +1)^{-k}\left( 1 - \Delta_v + |v|^2)^{k}\right)
\eeas
that $S_3 \in \vL(\vG^{k,s}, \vG^{k+2,s})$ for any $k \in \bZ$ and $s\in \bR$. A complex interpolation gives that $S_3 \in \vL(\vG^{r,s}, \vG^{r+2,s})$ for $r,s \in \bR $.
\ef


\begin{thebibliography}{99}

\bibitem{aw} M. Aafarani, X.P. Wang, {\it Gevrey estimates of the resolvent and sub-exponential time-decay for the heat and Schrödinger semigroups. II.}  J. Differential Equations \textbf{316} (2022), 387-424.

\bibitem{ag1} S. Agmon, {\it Spectral properties of Schrödinger operators and scattering theory}. Ann. Scuola Norm. Sup. Pisa Cl. Sci. (4) {\bf 2} (1975), no. 2, 151-218

\bibitem{ag2} S. Agmon,  {\it Lectures on Exponential Decay of Solutions to Second  Order Elliptic Equations}, Princeton Univ. Press, Princeton NJ, 1982.

\bibitem{mbs} M. Ben Said, {\it Kramers-Fokker-Planck operators with homogeneous potentials},   Math. Methods in the Appl. Sci.,  {\bf 45} (2) (2022) ,  914-927



\bibitem{d} E. B. Davies, {\it Linear Operators and Their Spectra}, Cambridge Univ. Press, London, 2007.



\bibitem{ff} J. Faupin, N. Frantz, {\it Spectral decomposition of some non-self-adjoint operators},
Annales Henri Lebesgue, {\bf 6 }(2023),  1115-1167.

\bibitem{ffro} J. Faupin and J. Fröhlich, {\it Asymptotic completeness in dissipative scattering theory}, Adv. Math., {\bf 340 }(2018), 300-362.

\bibitem{fh} R. Froese and  I. Herbst,
  {\it Exponential bounds and absence of positive eigenvalues for $N$-body {S}chr\"odinger operators},
 Commun. Math. Phys. {\bf 87}  (1982), 429--447.




\bibitem{hln} B. Helffer,  F. Nier, {\it  Hypoelliptic estimates and spectral theory for Fokker-Planck operators and Witten Laplacians}. Lecture Notes in Mathematics, 1862. Springer-Verlag, Berlin, 2005. x+209 pp. ISBN: 3-540-24200-7



\bibitem{hor} L.  Hörmander, {\it The analysis of linear partial differential operators. III. Pseudo-differential operators. }  Classics in Mathematics. Springer, Berlin, 2007.

\bibitem{ik} H. Isozaki, H. Kitada, {\it Micro-local resolvent estimates for 2-body Schr\"odinger operators}, J. Funct. Anal., {\bf 57} (1984), 270-300.

\bibitem{ky} T. Kako, K.  Yajima, {\it Spectral and scattering theory for a class of non-self-adjoint operators.}  Sci. Papers College Gen. Ed. Univ. Tokyo,  {\bf 26 } (1976), no. 2, 73-89



\bibitem{k} T. Kato, {\it Perturbation Theory for Linear Operators}, Springer, Berlin, 1984.

\bibitem{lz} T. Li, Z. Zhang, {\it Large time behaviour for the Fokker-Planck equation with general potential}, Sci. China, Mathematics, {\bf 61} (2018), Iss. 1, pp 137-150.

\bibitem{lwx}  W.X. Li, {\it Global hypoellipticity and compactness of resolvent for Fokker-Planck operator},  Ann. Sc. Norm. Super. Pisa Cl. Sci. (5) {\bf 11} (2012), no. 4, 789-815.

\bibitem{lwx2} W.X. Li, {\it Compactness criteria for the resolvent of the Fokker-Planck operator}, Ann. Sc. Norm. Super. Pisa Cl. Sci. (5) {\bf 18} (2018), no. 1, 119-143.



\bibitem{nw} R. Novak, X.P. Wang, {\it On the Kramers-Fokker-Planck equation with decreasing potentials in dimension one}, J. of Spectral Theory,
{\bf 10}(1) (2020),  1-32.


\bibitem{ol} F. W. J. Olver, {\it Asymptotics and Special Functions}, A. K. Peters Classics, Massachusetts, 1997.




\bibitem{risc} H. Risken,  {\it  The Fokker-Planck equation, Methods of solutions and applications}.  Springer, Berlin, 1989.

\bibitem{sa} Y. Saito, {\it The principle of limiting absorption for the nonself-adjoint Schrödinger operator in $R\sp{N}(N\neq 2)$.}  Publ. Res. Inst. Math. Sci.  {\bf 9 }  (1973/74), 397-428.



 \bibitem{sch} J. T.  Schwartz,  {\it Some non-self-adjoint operators} , Comm. Pure and Appl. Math., {\bf 13}(1960), 609-639.








\bibitem{w1} X. P. Wang, {\it Time-decay of semigroups generated by dissipative Schrödinger operators}. J. Differential Equations {\bf 253} (2012), no. 12, 3523-3542.


\bibitem{w2} X. P. Wang, {\it Large-time asymptotics of solutions to the Kramers-Fokker-Planck equation with a short-range potential.}  Comm. Math. Phys. {\bf 336 } (2015), no. 3, 1435-1471.

\bibitem{w3} X.P. Wang,  {\it Gevrey estimates of the resolvent and sub-exponential time-decay for the heat and Schrödinger semigroups}, J. Math. Pures Appl. \textbf{135} no. 9 (2020), 284-338

\bibitem{wz} X.P. Wang, L. Zhu, {\it Global-in-time $L^p-L^q$ estimaes for solutions of the Kramers-Fokker-Planck equation}, Commun. Math. Res., {\bf 38}(4) (2022), 560-578.
\end{thebibliography}
\end{document}